\documentclass[11pt,a4paper,reqno]{amsart}%reqno for right side numbering equation in RHs
\usepackage[utf8]{inputenc}
\usepackage[T1]{fontenc}
\usepackage{amsaddr}
\usepackage{amsmath}
\usepackage{xcolor}
\usepackage{soul}
\usepackage{graphicx}
\usepackage{wrapfig}
\usepackage{subcaption}
\usepackage[nopar]{lipsum}
\usepackage{amsthm}%for proof environment
\usepackage{amssymb}
\usepackage{euscript}
\usepackage{enumitem}%This package provides user control over the layout of the three basic list environments: enumerate, itemize and description.
\usepackage{mathtools}% for beautiful math.google for more info
\usepackage{soul}%for using \hl for highlighting
\usepackage{amscd}% for arrow diagram
\usepackage{float}
\usepackage[pagebackref,colorlinks=true,
linkcolor=blue,
urlcolor=red,colorlinks,
citecolor=green]{hyperref}
\usepackage{cleveref}
\usepackage{geometry}
\usepackage{pgf,tikz,pgfplots}
\usepackage{adjustbox}
\usepackage{setspace}
\usepackage[numbers,sort&compress]{natbib}

\allowdisplaybreaks[4]
\usetikzlibrary{calc}
\usetikzlibrary{patterns}
\pgfplotsset{compat=1.14}
\numberwithin{equation}{section}%If in an article you would like equations numbered by section
\newtheorem{thm}{Theorem}[section]
\newtheorem{prop}[thm]{Proposition}
\newtheorem{lem}[thm]{Lemma}

\newtheorem{exm}[thm]{Example}
\newtheorem{df}[thm]{Definition}

\newcommand {\R} {\mathbb{R}}
\newcommand {\n} {\mathbb{N}}

\newcommand {\dd} {\mathfrak{D}_d(H,g,f)}
\newcommand {\dc} {\mathfrak{D}_c(H,g,f)}
\makeatletter
\@namedef{subjclassname@1991}{Subject}
\@namedef{subjclassname@2000}{Subject}
\@namedef{subjclassname@2010}{Subject}
\@namedef{subjclassname@2020}{Subject}
\makeatother
\begin{document}
	\title[Multi-body interactions using Hypergraphs]{A study of diffusion in network with multi-body interactions using Hypergraphs}
	%\title[Multi-body interactions and synchronization in discrete dynamical networks ]{A study on Multi-body interactions and synchronization in discrete dynamical networks using spectra of hypergraphs}

	\email[ ]{\textit {{\scriptsize anirban.banerjee@iiserkol.ac.in, samironparui@gmail.com}}}
	
	\author[Banerjee, Parui]{Anirban Banerjee $\And$ Samiron Parui} 
	\address{Department of Mathematics and Statistics, Indian Institute of Science Education and Research Kolkata, Mohanpur-741246, India
		%\footnote{anirban.banerjee@iiserkol.ac.in, samironparui@gmail.com.}
	}
	
	\date{\today}
	
	
	\subjclass[2020]{Primary 05C82,  %Small world graphs, complex networks (graph-theoretic aspects)
		05C65 % Hypergraphs
		; Secondary 37C99, %Smooth dynamical systems 
		05C50, %Graphs and linear algebra (matrices, eigenvalues,etc.), 
		39A12, %Discrete version of topics in analysis
		34D06, %Synchronization of solutions to ordinary differential equations
		92B25 %Biological rhythms and synchronization
	}
	%\doublespacing
	\maketitle
	\begin{abstract}
		Being cognizant of the abundance of multi-body interactions in various complex systems, here we investigate a possible way to incorporate multi-body interactions in dynamical networks.  Adopting hypergraph as the underlying architecture aids our proposed dynamical network models to go beyond the traditional archetype of only pairwise interactions. We introduce some matrices associated with hypergraphs to incorporate multi-body frameworks in dynamic networks. We illustrate the fact that the approximation of multi-body interactions by pairwise binary interactions, i.e. considering graph as the underlying architecture of the corresponding dynamical network may lead to a wrong conclusion to the study. Here we use weighted hypergraphs to deal with the multi-body interactions of variable weights. We study the possibility of global and local synchronization in discrete and continuous-time dynamical networks. Some real-world numerical illustrations are included at the end to reinforce our theoretical results.
	\end{abstract}
	
	%	\tableofcontents
	
	\section{Introduction}
	The study of the  dynamical network deals with the evolution of individual dynamical systems on the vertices of the underlying graph. This study involves many fundamental concepts from non-linear dynamics and spectral graph theory. It has applications in some crucial multidisciplinary research areas involving Chemistry, Computer science, Physics, Mathematics, Biology, Social science, and Information science  \cite{social1, phy2,phy&bio, chem1, chem2, phy1, ma1, bio1, socialmedia, all,7277006}. The conventional graph topology perspective fails to incorporate multi-body interplay in the real world. A similar problem arises to represent group formation in scientific collaboration networks and social networking platforms like Facebook, WhatsApp, etc. \cite{petri2014homological,petri2018simplicial,abs}.
	In many real-world interacting systems, the interactions are not pairwise but involve a larger number of vertices at a time \cite{benson2018simplicial}. Therefore, in this article, we deal with the following two questions. (a) How can one incorporate multi-body interaction  in the mathematical models of the  dynamical networks? (b)Do the conditions for the occurrence of a specific phenomenon in dynamical networks determined from the standard paradigm of pairwise interaction also hold in the framework beyond pairwise interaction?

	The first question is addressed here by adopting the hypergraph as the underlying architecture of the dynamical network. The word phenomenon covers a large range in the second question. However, in this article, we restrict ourselves to the underlying hypergraphs' diffusive influence and resulting synchronisation in discrete time and continuous time cases. We also study the stability of the synchronisation in a dynamical network with multi-body interactions here.
	
	Though many studies have been done on synchronisations of the trajectories of a dynamical network with graph topology, only a few significant contributions have been made so far on the same for dynamical systems on hypergraphs. 
	
	Before further discussing synchronisation in dynamical networks with higher-order interactions, we would like to elucidate some terminology-related ambiguity.
	In \cite{multilayer}, authors have studied synchronisation on hyper-network, combining two or more graphs or multi-layer networks. For more details and references on multi-layer and multiplex networks, readers can see \cite{de2013mathematical, rakshit2020intralayer}. 
	In literature, sometimes, multi-layer networks and multiplex networks are called hyper-network \cite{multilayer}. The underlying structure of multi-layer networks and multiplex networks are graphs or combinations of graphs and  should not be confused with hypergraphs. 
	Throughout this article, we consider hypergraphs, a generalization of graphs, where each hyperedge can be any subset of the vertex set containing at least two elements.
	%, but our systems on hypergraphs are hypergraphs, which are entirely different from multi-layer networks. 
	In 2014, the first attempt was made to analyze synchronization in dynamical systems  on hypergraphs,  in which the authors have used continuous-time dynamical systems. It analyzed local synchronization with $3$-uniform hypergraphs \cite{first}.  A recent study has been made on local stability analysis of unweighted continuous-time dynamical systems on  hypergraphs in \cite{mulas2020coupled}. The diffusion matrix used here is different from ours.
	
	A \textit{weighted hypergraph} $H$ is the ordered triple $(V(H), E(H),w)$ such that $V(H)$ is a non-empty subset, $E(H)\subseteq P(V(H))$, where $P(V(H))$ is the power set of $V(H)$, and $w_H:E(H)\to (0,\infty)$ is a function. We refer to $V(H)$, $E(H)$, and $w$ as the \textit{set of vertices}, \textit{set of hyperedges}, and the \textit{weight of the hyperedges} of $G$, respectively.
	Each element of $V(H)$ and $E(H)$ is \textit{vertex} and \textit{hyperedge} in $H$, respectively. For $v\in V(H)$, \textit{the star of $v$} is $E_v(H)=:\{e\in E(H):v\in e\}$.
	A vertex measure on $H$ is a positive valued function $m_H:V(H)\to (0,\infty)$.
	A  hypergraph $H$ is connected if  for any two vertices $v_1, v_l\in V(H)$, there exists a sequence  of vertices $v_1, v_2,\cdots, v_l$, such that, $v_i,v_{i-1}\in e_i$, for some $e_i\in E(H)$, $i=1, \dots, l-1$. The rank, $rk(H)$, and co-rank, $cr(H)$ of a hypergraph $H$ are defined as $cr(H)=\min\{|e|:e\in E(H)\}$, and $rk(H)=\max\{|e|:e\in E(H)\}$.
	We call a hypergraph \textit{$m$-uniform hypergraph} if $|e|=m$, for all $e\in E(H)$. A $2$-uniform hypergraph is called a \textit{graph}. A hypergraph is represented by hypermatrix and tensors. The adjacency hypermatrix of a $m$-uniform hypergraph $H$ \cite{MR2900714, gen-spectra} is $\mathcal{A}_H=(a_{v_1,\ldots,v_m})_{v_i\in V(H);i=1,\ldots,m}$, where 
	$$a_{v_1,\ldots,v_m}=
	\begin{cases}
		\frac{w_H(e)}{(m-1)!}& \text{~if~}e=\{v_1,\ldots,v_m\}\in E(H)\\
		0&\text{~otherwise.}
	\end{cases}$$
	The action of a $m$-uniform hypergraph $H$ with vertex measure $m_H$, and hyperedge weight $w_H$ on $\mathbb{R}^{V(H)}$ is represented as a multi-linear function $\mathcal{A}_{(H,m_H,w_H)}:\mathbb{R}^{V(H)}\to\mathbb{R}^{V(H)} $, defined as
	$$(\mathcal{A}_{(H,m_H,w_H)}x)(v)=\sum\limits_{v_2,\ldots ,v_m\in V(H)}\frac{a_{v_2\ldots v_m}}{m_H(v)}x(v_2)\ldots x(v_m).$$
	A general (non-uniform) hypergraph $H$, the set of hyperedges, $E(H)=\bigcup\limits_{i=cr(H)}^{rk(H)}E^i(H)$, where $E^i(H)=\{e\in E(H):|e|=i\}$ and we refer the $i$-uniform hypergraph $H_i=(V(H),E^i(H))$ as the \textit{$i$-uniform layer} of $H$.
	For a general hypergraph $H$, with vertex measure $m_H$, and hyperedge weight $w_H$, the function $\mathcal{A}_{(H,m_H,w_H)}:\mathbb{R}^{V(H)}\to\mathbb{R}^{V(H)} $ is defined as $\mathcal{A}_{(H,m_H,w_H)}=\sum\limits_{i=cr(H)}^{rk(H)}\mathcal{A}_{(H_i,m_{H_i},w_{H_i})}$.
	
	%define hypergraph and graph here.
	
	\section{Multi-body interactions in  dynamical networks}
	
	A \textit{dynamical system} is an ordered triple $(X,\mathbb{T}, F)$, where the state space $X$ is a non-empty set, the domain of time $\mathbb{T}$ is a semigroup with identity $0_\mathbb{T}$, and $ F:\mathbb{T}\times X\to X$ is a function that describes the evolution of a parameter on the state space in such a way that if at time $0_{\mathbb{T}}$ the state of the parameter is $x_0\in X $ then at time $t\in \mathbb{T}$ the state would be $x_t=F(t,x_0)$. For $t_1,t_2\in\mathbb{T}$, at time $t_1+t_2$, the state of the dynamical system $x_{t_1+t_2}=F(t_2,(F(t_1,x_0)))$, where $+$ is the binary operation on the semigroup $\mathbb{T}$. The map $\gamma:\mathbb{T}\to X$, defined as $t\mapsto x_t$, is called the trajectories of the dynamical systems. We denote the trajectory as $ \{x_t\}_{t\in \mathbb{T}}$. Here, for our work,  we assume $X=\mathbb R$.
	For continuous-time, the domain $\mathbb{T}$ is the set of all non-negative reals, whereas, for  discrete-time dynamical systems, it is the set of all non-negative integers, i.e., $\mathbb{T}=\mathbb{N}\cup\{0\}$. Let $G$ be a graph, and there be a dynamical system $(X,\mathbb{T},F_v)$ on each vertex $v\in V(G)$ such that each edge $\{u,v\}\in E(G)$ acts as the coupling between $(X,\mathbb{T},F_u)$ and $(X,\mathbb{T},F_v)$. The collection of dynamical systems $\{(X,\mathbb{T},F_v)\}_{v\in V(G)}$ with all the edge-couplings together called a dynamical network on $G$.
	\subsection{Discrete dynamical network with multi-body interactions}
	Let $G$ be a \textit{weighted graph}, where $V(G)$ and $E(G)$ are the set of \textit{vertices} and \textit{edges}, respectively. The weight of $G$ is a positive valued function $w_G:E(G)\to(0,\infty)$. Suppose that on each vertex $v\in V(H)$, there is a discrete dynamical system $(\mathbb{R},\mathbb{N}\cup \{0\},F_v)$ described by the iteration equation $x_{t+1}(v)=f(x_{t}(v))$, where $f:\mathbb{R}\to\mathbb{R}$. Thus, $F_v(t,x_0(v))=f^t(x_0(v))$. Suppose that 
	%Since the same function $f$ describes discrete-dynamical system on each vertices, the function 
	$f_G:\mathbb{R}^{V(G)}\to \mathbb{R}^{V(G)}$, defined by $(f_G(x_t))(v)=(f(x_t(v)))$, describes the collection of all the dynamical systems $\{(\mathbb{R},\mathbb{N}\cup \{0\},F_v)\}_{v\in V(G)}$, where $\mathbb{R}^{V(G)}$ is the set of all functions from $V(G)$ to $\mathbb{R}$. 
	Since each $\{u,v\}\in E(G)$ represents a coupling between the dynamical systems on the vertices $u$ and $v$, the dynamical network on $G$ is given by the following equation.
	\begin{align}\label{graph}
		x_{t+1}(v)=g(x_t(v)) +\epsilon \sum\limits_{ u\in V(G)}A_{uv}\left(f(x_t(u))-f(x_t(v))\right),
	\end{align} 
	where  $f:\mathbb{R}\to \mathbb{R}$ and  $g:\mathbb{R}\to \mathbb{R}$ are the functions describing the dynamical systems on each vertex, and $\epsilon$ is the coupling strength of the dynamical networks. For all $u,v\in V(G)$, $A_{uv}=\frac{w_G(\{u,v\})}{m_G(u)}$, if $\{u,v\}\in E(H)$, otherwise, $A_{uv}=0$, where $m_G:V(H)\to (0,\infty)$ is a discrete measure on the vertices of $G$. We can also express \Cref{graph} as follows.
	\begin{align}\label{graph-vector}
		x_{t+1}=g_G(x_t)+\epsilon L_G(f_G(x_t)),
	\end{align}
	where $L_G:\mathbb{R}^{V(H)}\to \mathbb{R}^{V(H)}$ and is defined as $(L_G(x))(v)=\sum\limits_{u\in V(G)}A_{uv}\left(f(x_t(u))-f(x_t(v))\right)$. Here, on each vertex, the dynamical system $(X,\mathbb{T},F_v)$ has two components. One is an \textit{interacting component} controlled by the function $f$, that affects the other dynamical systems, which are connected with $(X,\mathbb{T}, F_v)$ by hyperedge couplings in $E_v(G)$. We refer the other component of $(X,\mathbb{T},F_v)$ as a \textit{non-interacting component}, which is controlled by the function $g$ and does not affect any other dynamical system in the dynamical network. Since the discrete dynamical network given by \Cref{graph} or by \Cref{graph-vector} depends only on the functions $f,g$, and the graph $G$, we denote it by $\mathfrak{D}_d(G,g,f)$.
	Similar dynamical network models have been reported in the literature \cite{1,li2004synchronization,rangarajan2002stability,lu2008synchronization} in the last few decades. In \cite{1}, the vertex-measure $m_G$ is considered as $m_G(v)=|E_v(G)|$, whereas in \cite{li2004synchronization}, $m_G(v)=1$ for all $v\in V(G)$. Some models do not distinguish the interacting and non-interacting components and assume $g=f$ \cite{1}. 
	
	Since the discrete dynamical network model mentioned above is based on a graph, thus can not incorporate beyond binary interactions. So to incorporate multi-body interactions, we use a hypergraph $H$ as the underlying architecture of the discrete dynamical network. Suppose that $\{(X,\mathbb{T}, F_v)\}_{v\in V(H)}$ is a collection of dynamical systems such that $\{v_1,\ldots,v_m\}\in E(H)$ if and only if there is a multi-body interaction among the collection of dynamical systems $\{(X,\mathbb{T}, F_{v_i})\}_{i=1}^m$, i.e., here, each hyperedge acts as a multi-body coupling among the dynamical systems. Thus, the hypergraph $H$ represents a dynamical network that allows multi-body interactions, and we denote the corresponding dynamical network as $\dd$. Now the obvious question is- what is the model for this dynamical network on $H$? Which linear operator substitutes $L_G$ when the graph $G$ is replaced by a hypergraph $H$ in the architecture of the dynamical network? We try to find the answer in the following Theorem.
	For any hypergraph $H$, the \emph{$i$-uniform layer of} $H$ is a hypergraph $H_i$ such that $V(H_i)=V(H)$, and $E(H_i)=\{e\in E(H):|e|=i\}$.
	\begin{thm}\label{hypergraph-discrete}
		Let $H$ be a weighted hypergraph with hyperedge weight $w_H$ and vertex measure $m_H$, and $\dd$ be a discrete dynamical network. For any $t\in \mathbb{T}=\n\cup\{0\}$,
		\begin{align}\label{hypergraph-vector}
			x_{t+1}=g_H(x_t)+\epsilon C_{(H,m_H,w_H)}(f_H(x_t)),
		\end{align}
		where $$C_{(H,m_H,w_H)}:\mathbb{R}^{V(H)}\to \mathbb{R}^{V(H)}$$ is a linear operator defined by $C_{(H,m_H,w_H)}=\sum\limits_{i=cr(H)}^{rk(H)}B^{(i)}_{(H_i,m_{H_i},w_{H_i})}$, where $H_i$ is {$i$-uniform layer of $H$}, and $(B^{(i)}_{(H_i,m_{H_i},w_{H_i})}(x))(v)=\frac{1}{i}\sum\limits_{e\in E_v(H_i)}\frac{w_{H_i}(e)}{m_{H_i}(v)}\sum\limits_{u\in e}(x(u)-x(v))$.
		
	\end{thm}
	\begin{proof}
		First consider the case when the underlying hypergraph is an $m$-uniform hypergraph.
		To develop the generalized notion for a discrete dynamical network with a hypergraph as its underlying topology,  the binary difference term $\left(f(x_t(u))-f(x_t(v))\right) $ in \Cref{graph} needs to be replaced by a multi-nary influence term corresponding to a hyperedge. For the multi-nary influence of the hyperedge-coupling corresponding to a hyperedge \sloppy $e=\{u,u_2,\ldots, u_m\}$ on the dynamical system on the vertex $u$ we consider the term $\sum\limits_{ u_2,\ldots, u_m\in V}\frac{a_{uu_2\ldots u_m}}{m_H(u)}\left(\frac{\sum\limits_{j=2}^m {f}(x_t(u_j))}{m-1}-{f}(x_t(u))\right)$. Thus the  discrete dynamical network  with $m$-uniform multi-nary interaction can be expressed as, 
		\begin{align}\label{m-hypergraph}
			x_{t+1}(v)=g(x_t(v)) +\epsilon \sum\limits_{ u_2,\ldots, u_m\in V}\frac{a_{uu_2\ldots u_m}}{m_H(u)}\left(\frac{1}{m-1}\sum\limits_{j=2}^m {f}(x_t(u_j))-{f}(x_t(u))\right).
		\end{align}
		To express the above equation more concisely, we need an operator in place of $L_G$ in \Cref{graph-vector}.
		Since $\frac{1}{m-1}\sum\limits_{j=2}^m {f}(x_t(u_j))-{f}(x_t(u))=\frac{m}{m-1}(\frac{1}{m}(\sum\limits_{j=2}^m {f}(x_t(u_j))+f(x_t(u)))-{f}(x_t(u))$, contribution of each $e\in E_u(H)$ to the sum $\sum\limits_{ u_2,\ldots, u_m\in V}a_{uu_2\ldots u_m}\left(\frac{1}{m-1}\sum\limits_{j=2}^m {f}(x_t(u_j))-{f}(x_t(u))\right)$ is $w_H(e)\frac{m}{m-1}(\frac{1}{m}(\sum\limits_{v\in e}f(x_t(v)))-{f}(x_t(u))) $. Therefore, $\sum\limits_{ u_2,\ldots, u_m\in V}\frac{a_{uu_2\ldots u_m}}{m_H(u)}\left(\frac{1}{m-1}\sum\limits_{j=2}^m {f}(x_t(u_j))-{f}(x_t(u))\right)= \frac{m}{(m-1)}\sum\limits_{e\in E_u(H)}\frac{w_H(e)}{m_H(u)}(\frac{1}{m}(\sum\limits_{v\in e}f(x_t(v)))-{f}(x_t(u))) $.
		
		For the functions $sum_V:\mathbb{R}^{V(H)}\to \mathbb{R}^{E(H)} $, and $sum_E:\mathbb{R}^{E(H)}\to \mathbb{R}^{V(H)} $ defined as $(sum_V(x))(e) =\sum\limits_{v\in e}x(v)$, and $(sum_E(\alpha))(v) =\frac{1}{m_H(v)}\sum\limits_{e\in E_v(H)}w_H(e)\alpha(e)$, respectively, the function $Q^{(m)}_{(H,m_H,w_H)}=\frac{1}{m}(sum_E\circ sum_V):\mathbb{R}^{V(H)}\to \mathbb{R}^{V(H)} $ is given by 
		$$(Q^{(m)}_{(H,m_H,w_H)}(x))(v)=\frac{1}{m}\sum\limits_{e\in E_v(H)}\frac{w_H(e)}{m_H(v)}\sum\limits_{u\in e}x(u)$$
		for all $x\in \mathbb{R}^{V(H)} $, $v\in V(H)$. 
		We define $B^{(m)}_{(H,m_H,w_H)}:\mathbb{R}^{V(H)}\to \mathbb{R}^{V(H)} $ as $(B^{(m)}_{(H,m_H,w_H)}(x))(v)=(Q^{(m)}_{(H,m_H,w_H)}(x))(v)-\sum\limits_{ e\in E_v(H)}\frac{w_H(e)}{m_H(v)}x(v)$. Thus,
		$$(B^{(m)}_{(H,m_H,w_H)}(x))(v)=\frac{1}{m}\sum\limits_{e\in E_v(H)}\frac{w_H(e)}{m_H(v)}\sum\limits_{u\in e}(x(u)-x(v)),$$
		and therefore, \Cref{m-hypergraph} becomes $	x_{t+1}(v)=g(x_t(v)) +\epsilon \frac{m}{m-1}(B^{(m)}_{(H,m_H,w_H)}(f_H(x_t)))(v)$. Thus, the discrete dynamical network with an $m$-uniform hypergraph as its underlying architecture is 
		\begin{align}\label{m-hypergraph-B}
			x_{t+1}=g_H(x_t) +\epsilon \frac{m}{m-1}B^{(m)}_{(H,m_H,w_H)}(f_H(x_t)).
		\end{align}
		Now we are in a position to consider the case when the underlying hypergraph $H$ of the dynamical network is non-uniform. If $H$ is non-uniform then the set of hyperedges, $E(H)=\bigcup\limits_{i=cr(H)}^{rk(H)}E^i(H)$, where $E^i(H)=\{e\in E(H):|e|=i\}$.
		For a general hypergraph $H$, with vertex measure $m_H$, and hyperedge weight $w_H$, we define the function $C_{(H,m_H,w_H)}:\mathbb{R}^{V(H)}\to\mathbb{R}^{V(H)} $ as $C_{(H,m_H,w_H)}=\sum\limits_{i=cr(H)}^{rk(H)}\frac{i}{i-1}B^{(i)}_{(H_i,m_{H_i},w_{H_i})}$, and the discrete dynamical network model is
		\begin{align}\label{m-hypergraph-c}
			x_{t+1}=g_H(x_t) +\epsilon C_{(H,m_H,w_H)}(f_H(x_t)).
		\end{align}
	\end{proof}
	%%%%%%%%%%%%%%%%%%%%%%%%%%%%%%%%%
	%%%%%%%%%%%%%%%%%%%%%%%%%%%%%%%%%
	In the \Cref{hypergraph-discrete}, we have incorporated the multi-body interaction in the discrete dynamical network using the operator $C_{(H.m_H,w_H)}$ associated with the underlying hypergraph. We refer the operator $C_{(H.m_H,w_H)}$ as\textit{ multi-body interaction operator}.  %Since the hyperedges act as the couplings in the discrete dynamical network, {\color{red} the multi-body interaction should be the combined action of the hyperedges.} Let $H$ be the underlying hypergraph of the discrete dynamical network. 
	For each $e\in E(H)$, we incorporate the action of $e$ using the function $L_e:\mathbb{R}^{V(H)}\to \mathbb{R}^{V(H)} $, and is defined as 
	$$(L_e(x))(v)=
	\begin{cases}
		\sum\limits_{u\in e}\frac{1}{m_H(u)}\frac{1}{|e|-1}(x(u)-x(v)) &\text{~if~} v\in e\\
		0&\text{~otherwise,}
	\end{cases}
	$$
	for all $x\in \mathbb{R}^{V(H)} $, $v\in V(H)$.
	We encode the combined action of all the hyperedges using the operator $ L_{(H,m_H,w_H)}:\mathbb{R}^{V(H)}\to \mathbb{R}^{V(H)} $ defined by $L_{(H,m_H,w_H)}=\sum\limits_{e\in E(H)}w_H(e)L_e$. In the following result, we show that the action of the multi-body interaction operator is equal to the combined action of the collection of operators $\{L_e\}_{e\in E(H)}$.
	\begin{thm}
		Let $\dd$ be a discrete dynamical network  with hyperedge weight $w_H$ and vertex measure $m_H$ of the underlying hypergraph $H$. The multi-body interaction operator, $C_{(H,m_H,w_H)}=L_{(H,m_H,w_H)}$.
	\end{thm}
	\begin{proof}
		Since $C_{(H,m_H,w_H)}=\sum\limits_{i=cr(H)}^{rk(H)}\frac{|i|}{|i|-1}B^{(m)}_{(H_i,m_{H_i},w_{H_i})}$, where $$(B^{(i)}_{(H_i,m_{H_i},w_{H_i})}(x))(v)=\frac{1}{i}\sum\limits_{e\in E_v(H_i)}\frac{w_{H_i}(e)}{m_{H_i}(v)}\sum\limits_{u\in e}(x(u)-x(v))$$   for all $x\in \mathbb{R}^{V(H)} $, and $v\in V(H)$,
		\begin{align*}
			( C_{(H,m_H,w_H)}(x))(v)&=\sum\limits_{i=cr(H)}^{rk(H)}\frac{1}{i-1}\sum\limits_{e\in E_v(H_i)}\frac{w_{H_i}(e)}{m_{H_i}(v)}\sum\limits_{u\in e}(x(u)-x(v))\\
			&=\sum\limits_{e\in E_v(H)}\frac{w_{H}(e)}{m_{H}(v)}\frac{1}{|e|-1}\sum\limits_{u\in e}(x(u)-x(v))\\
			&=\sum\limits_{e\in E(H)}w_H(e)(L_e(x))(v)=(L_{(H,m_H,w_H)}(x))(v).
		\end{align*}  
		This completes the proof.
	\end{proof}
	\begin{df}
		We call an operator $A:\mathbb{R}^{V(H)}\to \mathbb{R}^{V(H)} $ a diffusion operator if $A$ has the following property.
		\begin{enumerate}
			\item $0$ is an eigenvalue of $A$ and the corresponding eigenvector is $\chi_{V(H)}$.
			\item Other than $0$ all the eigenvalues of $A$ are negative.
		\end{enumerate}
	\end{df}
	If $x(t)\in \mathbb{R}^{V(H)}$ is a solution of the differential equation $\dot{x}=Ax$, where $A$ is a diffusion operator, then $\lim\limits_{t\to\infty}x(t)=c\chi_{V(H)}$ for some $c\in\R$. thus under the action of the operator $A$, a diffusion process is taking place, and all the components of $x(t)$ tend to be equal as $t\to \infty$. Moreover, as we know, the diffusion process ends at density equality; here, we also have $A(c\chi_{V(H)})=0$. 
	For any hypergraph $H$, with hyperedge weight $w_H$ and vertex measure $m_H$, there exists two inner products $(\cdot,\cdot)_V$, and $(\cdot,\cdot)_E$ defined as $(x,y)_V=\sum\limits_{v\in V(H)}m_H(v)x(v)y(v)$, and $(\alpha,\beta)_E=\sum\limits_{e\in E(H)}w_H(e)\alpha(e)\beta(e)$, where $x,y\in \mathbb{R}^{V(H)}$, $\alpha,\beta\in \mathbb{R}^{E(H)}$.%, and $\langle\cdot,\cdot\rangle$ is the usual inner product in $\mathbb{R}$. 
	\begin{lem}\label{lem-selfadj}
		For any hypergraph $H$, with hyperedge weight $w_H$, and vertex measure $m_H$, the adjoint of $sum_V$ is $sum_E$.
	\end{lem}
	\begin{proof}
		For $x\in \mathbb{R}^{V(H)}$, $\alpha\in \mathbb{R}^{E(H)}$, 
		\begin{align*}
			(sum_V(x),\alpha)_V &=\sum\limits_{e\in E(H)}w_H(e)(sum_V(x))(e)\alpha(e) \\
			&=\sum\limits_{e\in E(H)}w_H(e) \sum\limits_{v\in e}x(v)\alpha(e) \\
%			&=\sum\limits_{e\in E(H)}w_H(e) \sum\limits_{v\in e}\langle x(v),\alpha(e) \rangle\\
			&=\sum\limits_{v\in V(H)}m_H(v) x(v)\sum\limits_{e\in E_v(H)}\frac{w_H(e)}{m_H(v)}\alpha(e)
			\\
			&=(x,sum_E(\alpha))_E.
		\end{align*}
		This completes the proof.
	\end{proof}
	\begin{lem}\label{lem-nd}
		Let $H$ be a hypergraph with hyperedge weight $w_H$ and vertex measure $m_H$. For any $x\in \mathbb{R}^{V(H)}$,
		$$((C_{(H,m_H,w_H)}x),x)_V=-\sum\limits_{v\in V(H)}\frac{w_H(e)}{2(|e|-1)}\sum\limits_{u,v\in e}|x(v)-x(u)|^2.$$
	\end{lem}
	\begin{proof}
		For any $x\in \mathbb{R}^{V(H)}$,
		\begin{align*}
			((C_{(H,m_H,w_H)}x),x)_V&=\sum\limits_{v\in V(H)}m_H(v) (C_{(H,m_H,w_H)}x)(v)) x(v)\\
			&=\sum\limits_{v\in V(H)}\sum\limits_{e\in E_v(H)}\frac{w_H(e)}{|e|-1}\sum\limits_{u\in e} (x(u)-x(v))x(v)\\
			&=-\sum\limits_{e\in E(H)}\frac{w_H(e)}{2(|e|-1)}\sum\limits_{u,v\in e} (x(v)-x(u))^2.
		\end{align*}
		This completes the proof.  
	\end{proof}
	\begin{thm}
		For any connected hypergraph $H$, with hyperedge weight $w_H$ and vertex measure $m_H$, the operator $ C_{(H,m_H,w_H)}$ is a diffusion operator. 
	\end{thm}
	\begin{proof}
		For any $m$-uniform hypergraph, by \Cref{lem-selfadj}, being a composition of a linear map and its adjoint, the operator $Q^{(m)}_{(H,m_H,w_H)}$ is self-adjoint.
		The map $D^{(m)}_H:\mathbb{R}^{V(H)}\to \mathbb{R}^{V(H)}$ defined by $(D^{(m)}_H(x))(v)=\sum\limits_{ e\in E_v(H)}\frac{w_H(e)}{m_H(v)}x(v) $ is also self-adjoint. Thus $ B^{(m)}_{(H,m_H,w_H)}=Q^{(m)}_{(H,m_H,w_H)}-D^{(m)}_H$ is also self-adjoint. Therefore, for any non-uniform hypergraph $H$, $C_{(H,m_H,w_H)}=\sum\limits_{i=cr(H)}^{rk(H)}\frac{|i|}{|i|-1}B^{(m)}_{(H_i,m_{H_i},w_{H_i})}$ is self-adjoint.
		By \Cref{lem-nd}, $C_{(H,m_H,w_H)}$ is negative semidefinite and thus other than $0$, all the eigenvalues are negative. If $0$ is an eigenvalue of $C_{(H,m_H,w_H)}$ with eigenvector $x_0$, then by \Cref{lem-nd}, for a connected hypergraph $H$, $x_0(u)=x_0(v)$ for all $u,v\in V(H)$. Thus, $x_0=c\chi_{V(H)}$, for some $c\in\mathbb{R}$. This completes the proof.
	\end{proof}
	
	Note that for a hypergraph without a loop, the operator $C_{(H,m_H,w_H)}$ studied here becomes the general diffusion operator $(\mathcal{L}_{(H,\delta_{V(H)},\delta_{E(H)})}$ associated with a hypergraph $H$, introduced in \cite{up2021} for $m_H=\delta_{V(H)}$, and $w_H(e)=\delta_{E(H)}(e)\frac{|e|-1}{|e|^2}$ for all $e\in E(H)$.
	\subsection{Continuous dynamical network with multi-body interactions}
	Let $G$ be a graph with edge weight $w_G:E(G)\to (0,\infty)$ and vertex measure $m_H:V(G)\to (0,\infty)$. Suppose that there are identical continuous time dynamical system $(\mathbb{R},[0,\infty),F_v)$ on $v\in V(H)$ defined by the differential equation $\dot{x}_t(v)=f(x_t(v))$, where $f:\mathbb{R}\to\mathbb{R}$ is a differentiable function. That is $F_v(t,x_0)=s_{x_0}(t)$, where $s_{x_0}$ is the solution of the differential equation $\dot{x}_t(v)=f(x_t(v))$ with $ x_t=x_0$ when $t=0$. Thus, the function $f_G:\mathbb{R}^{V(H)}\to\mathbb{R}^{V(H)}$ defined by $(f_G(x_t))(v)=(f(x_t(v)))$, describes the collection of all the dynamical systems $\{(\mathbb{R},[0,\infty),F_v)\}_{v\in V(G)}$.  Since each $\{u,v\}\in E(G)$ works as a coupling between the dynamical systems on the vertices on $u$ and $v$, the dynamical network is given by the following equation.
	\begin{equation}\label{contn-graph}
		\Dot{x}_t=g_G(x_t)+\epsilon L_G( f_G(x_t)),
	\end{equation}
	where $f:\mathbb{R}\to \mathbb{R}$ and $g:\mathbb{R}\to \mathbb{R}$ are differentiable functions describing the dynamical systems on each vertex. Similar continuous time dynamical system models can be found in \cite{MR2178380,MR2265856} and references therein. Like the discrete case, if we want to replace the underlying graph $G$ by a hypergraph $H$ with hyperedge weight $w_H$ and the vertex measure $m_H$, then we have to replace the operator $L_G$ by the operator $C_{(H,m_H,w_H)}$. Thus, the continuous-time dynamical network model is given by the following equation.
	\begin{equation}\label{contn-w}
		\Dot{x}_t=g_H(x_t)+\epsilon C_{(H,m_H,w_H)}( f_H(x_t)).
	\end{equation}
	We denote this continuous dynamical network as $\dc$. 
	\section{Synchronization} Suppose $\dd$ (or $\dc$) is a discrete (or continuous) dynamical network.
	If $x_t(u)=x_t(v)$ for $u,v\in V(H)$ at time $t$ then we say $\dd$ (or $\dc$) \textit{synchronizes at time $t$}. Moreover, if $\lim\limits_{t\to\infty}|x_t(u)-x_t(v)|=0$, we say that the dynamical network \textit{synchronizes asymptotically}. Since we have constructed  linear operators associated with hypergraphs to incorporate the diffusion process in dynamical networks, we can now apply traditional techniques for the conventional dynamical network on graphs to study the synchronization in dynamical networks with multi-body interactions. 
	\subsection{Stability analysis of synchronization in discrete dynamical network} Let $\dd$ be a discrete dynamical network.
	The stability analysis of $\dd$ is the study of the question if any synchronized trajectory $\{x_t\}_{t\in\n\cup \{0\}}$ is stable under small perturbation $\eta_t\in{\R}^{V(H)} $. That is, for some $x_t,n_t\in {\R}^{V(H)}$ with  $x_t(v)=c_t(\in\mathbb{R})$, and $|\eta_0(v)|$ is very small for all $v\in V(H)$, and for all $t\in\n\cup\{0\}$, whether the trajectory $\{y_t\}_{t\in\n\cup \{0\}}$ achieve synchronization, where $y_t=x_t+\eta_t$ for all $t$. 
	
	Since $ C_{(H,m_H,w_H)}$ is self-adjoint, all the eigenvectors of the operator form an orthonormal basis of $\mathbb{R}^{V(H)}$. Let $\mathfrak{B}=\{z_i\}_{i=1}^{|V(H)|}$ be that orthonormal basis of $\mathbb{R}^{V(H)}$, which consists of the eigenfunction of $ C_{(H,m_H,w_H)}$. Thus, $x_t=\sum\limits_{z_i\in \mathfrak{B}}(x_t,z_i)_Vz_i$ for all $t$.
	\begin{prop}
		Let $\dd$ be a discrete dynamical network with the hyperedge weight and the vertex measure of $H$ are $w_H$, $m_H$, respectively, and $f,g$ are differentiable functions with bounded derivatives,  $\sup(f^\prime)=\sup\limits_{x\in\mathbb{R}}(\frac{df(x)}{dx})$, $\sup(g^\prime)=\sup\limits_{x\in\mathbb{R}}(\frac{dg(x)}{dx})$. If absolute values of all the eigenvalues of $ C_{(H,m_H,w_H)}$ belongs to the interval $(-\frac{1+|\sup(f^\prime)|}{\epsilon |\sup(g^\prime)|}, \frac{1-|\sup(f^\prime)|}{\epsilon |\sup(g^\prime)|})$ then synchronization in $\dd$ is stable under small perturbations.
	\end{prop}
	\begin{proof}

		Since $y_t=x_t+\eta_t$ for all $t\in\n\cup \{0\}$, by \Cref{hypergraph-vector} we have $\eta_{t+1}=g_H(y_t)-g_H(x_t)+\epsilon C_{(H,m_H,w_H)}(f_H(y_t)-f_H(x_t))$. Suppose that $\eta_{t+1}=(\eta_{t+1},z_i)_V$, thus, $\eta_{t+1}=\sum\limits_{z_i\in \mathfrak{B}}\eta^i_{t+1}z_i$.
		Now, by mean value theorem, $|(g_H(y_t)-g_H(x_t),z_i)_V|\le |\sup(g^\prime)||\eta_t^i|$. Therefore,
		\begin{align*}
			|\eta^i_{t+1}|&=|(g_H(y_t)-g_H(x_t),z_i)_V+\epsilon (C_{(H,m_H,w_H)}(f_H(y_t)-f_H(x_t)),z_i)_V|\\
			&\le (|\sup(f^\prime)|+\epsilon|\lambda_i||\sup(g^\prime)|)|\eta_t^i|.
		\end{align*}
		Thus, if $-1<(|\sup(f^\prime)|+\epsilon|\lambda_i||\sup(g^\prime)|)
		<1$, then $ \lim\limits_{t\to\infty}\eta_t^i=0$. Thus, the result follows.
	\end{proof}
	If the perturbation $\eta_t$ is very small then $f(y_t(v))=f(x_t(v))+f^\prime (c_t)\eta_t(v)$ for all $v\in V(H)$. This leads us to the following result.
	\begin{prop}
		Let $\dd$ be a discrete dynamical network with the hyperedge weight and the vertex measure of $H$ are $w_H$, $m_H$, respectively, and $f,g$ are differentiable functions. If $\lim\limits_{t\to\infty}\sigma_i(t)<0$ for all $i=1,\ldots,|V(H)|$ then any synchronization in $\dd$ is stable under small perturbation, where $\sigma_i(t)=\frac{1}{t}\sum\limits_{i=1}^tlog_e|(f^\prime (c_t)+\epsilon\lambda_ig^\prime (c_t))|$.
	\end{prop}
	\begin{proof}
		Since the perturbations are small, for all $v\in V(H)$, $f(y_t(v))=f(x_t(v))+f^\prime (c_t)\eta_t(v)$. Thus, $(f_H(y_t),z_i)_V=(f_H(x_t),z_i)_V+f^\prime (c_t)\dot{\eta}_t^i$. Therefore, $\eta^i_{t+1}=(g_H(y_t)-g_H(x_t),z_i)_V+\epsilon (C_{(H,m_H,w_H)}(f_H(y_t)-f_H(x_t)),z_i)_V=(f^\prime (c_t)+\epsilon\lambda_ig^\prime (c_t))\eta_t^i$, and $|\eta^i_{t+1}
		|=\prod\limits_{i=0}^t|(f^\prime (c_t)+\epsilon\lambda_ig^\prime (c_t))|\eta_0^i=e^{t\sigma_i(t)}\eta_0^i$, where $\sigma_i(t)=\frac{1}{t}\sum\limits_{i=1}^tlog_e|(f^\prime (c_t)+\epsilon\lambda_ig^\prime (c_t))|$. Thus, if $\lim\limits_{t\to\infty}\sigma_i(t)<0$ then $\lim\limits_{t\to\infty}\eta^i_t=0$ and the result follows. 
	\end{proof}
	The above result has an interesting consequence when $f=g$.
	\begin{prop}
		Let $\mathfrak{D}_d(H,f,f)$ be a discrete dynamical network with the hyperedge weight and the vertex measure of $H$ are $w_H$, $m_H$, respectively, and $f$ is a differentiable function, and  $\sigma(t)=\frac{1}{t}\sum\limits_{i=1}^t\log_e|f^\prime(c_t)|$. If all the eigenvalues of $ C_{(H,m_H,w_H)} $ belongs to the interval $(-\frac{1}{\epsilon}(\frac{1}{e^{\sigma_\infty}}+1),\frac{1}{\epsilon}(\frac{1}{e^{\sigma_\infty}}-1))$ then any synchronization in  $\mathfrak{D}_d(H,f,f)$ is stable under small perturbation, where $\sigma_\infty=\lim\limits_{t\to\infty}\sigma(t)$.
	\end{prop}
	\begin{proof}
		Suppose that $\sigma(t)=\frac{1}{t}\sum\limits_{i=1}^t\log_e|f^\prime(c_t)|$.     Since for $f=g$,  $\sigma_i(t)=\sigma(t)\log_e|
		(1+\epsilon \lambda_i)|$, the condition of stability of synchronization in $\mathfrak{D}_d(H,f,f)$ is $e^{\sigma_{\infty}}|
		(1+\epsilon \lambda_i)|<1$, where $\sigma_\infty=\lim\limits_{t\to\infty}\sigma(t)$. Since $e^{\sigma_{\infty}}|
		(1+\epsilon \lambda_i)|<1$ is equivalent to $\lambda_i\in (-\frac{1}{\epsilon}(\frac{1}{e^{\sigma_0}}+1),\frac{1}{\epsilon}(\frac{1}{e^{\sigma_0}}-1))$, the result follows.
	\end{proof}
	\subsection{Global analysis of synchronization in discrete dynamical networks}
	Now we study the conditions that compelled the trajectories of $\dd$ to attain synchronization. Let $\{x_t\}_{t\in\n\cup\{0\}}$ be a trajectory of $\dd$. Let
	$ 	s_0=\frac{1}{|V(H)|}\sum_{u\in V(H)} x_0(u)$, and  $\{y_t\}_{t\in\n\cup\{0\}}$ be a trajectory of $\dd$ such that $y_0(v)=s_0$ for all $v\in V(H)$. Thus, $\{y_t\}_{t\in\n\cup\{0\}}$ is a synchronized trajectory with $y_{t+1}(v)=g(c_{t})$ for all $v\in V(H)$. If $\lim\limits_{t\to\infty}|x_t-y_t|=0$, then $\{x_t\}_{t\in\n\cup\{0\}}$ attain synchronization asymptotically. Thus, we have the following Proposition.
	\begin{prop}\label{lpshz-dis}
		Let $\dd$ be a discrete dynamical network. If $f$ and $g$ are Lipschitz functions with Lipschitz constant $k_f$ and $k_g$, respectively, and $[k_g+\epsilon \| C_{(H,m_H,w_H)}\|k_f]< 1, $ where $\| C_{(H,m_H,w_H)}\| $ is the operator norm of $ C_{(H,m_H,w_H)}$,
		then any trajectory of $\dd$ achieves synchronization asymptotically. Moreover, if $f=g$ and then the condition for the synchronization is $\|[I+\epsilon  C_{(H,m_H,w_H)}]\|< \frac{1}{k_f}$, where $I$ is the identity operator.
	\end{prop}
	\begin{proof}
		Since $|x_{t+1}-y_{t+1}|\le \|[k_g+\epsilon \| C_{(H,m_H,w_H)}\|k_f]|x_{t}-y_{t}|$, the result follows.
	\end{proof}
	The operator norm of a self-adjoint operator is the  maximum of the spectral radius of the operator. Thus, if $\lambda_{\max}$ is the eigenvalue of $C_{(H,m_H,w_H)}$ such that $[k_g+\epsilon \| \lambda_{\max}\|k_f]=\max\{[k_g+\epsilon \| \lambda\|k_f]:\lambda \in S(C_{(H,m_H,w_H)} )\}$, where $S(C_{(H,m_H,w_H)}$ is set of all the eigenvalues of $C_{(H,m_H,w_H)}$, then the condition for synchronization becomes $[k_g+\epsilon \| \lambda_{\max}\|k_f]<1$.
	\subsection{Stability analysis of synchronization in continuous dynamical network}
	Let $\{x_t\}_{t\in[0.\infty)}$ be a synchronized trajectories of $\dc$. Thus, for all $t$, $x_t=c_t\chi_{V(H)}$ for some $c_t\in \mathbb{R}$. If the synchronized trajectory is perturbed by the small initial perturbation $\eta_0$, then it becomes the perturbed trajectory $\{y_t\}_{t\in[0.\infty)}$, where $y_t=x_t+\eta_t$ for all $t\in [0,\infty)$. Thus, by \Cref{contn-w},
	$\Dot{\eta}_t=f_G(y_t)-f_G(x_t)+\epsilon C_{(H,m_H,w_H)}( g_H(y_t)- g_H(x_t))$. Suppose that $\eta
	_t^i=(\eta_t,z_i)_V$. Thus, $\dot{\eta}
	_t^i=({\dot\eta}_t,z_i)_V$. Since the perturbation is small, $g(y_t(v))=g(x_t(v))+g^\prime(c_t)\eta_t(v)$. Therefore, $( g_H(y_t)- g_H(x_t),z_i)_V=g^\prime(c_t)\eta_t^i$. Proceeding similarly we get $\dot{\eta}_t^i=({\dot\eta}_t,z_i)_V=g^\prime(c_t)\eta_t^i+\epsilon \lambda_if^\prime(c_t)\eta_t^i$. Therefore,
	$\eta_t^i=\eta_0^ie^{\int_{0}^{t} [g^\prime(c_r)+\epsilon \lambda_if^\prime(c_r)]\,dr }\le \eta_0^ie^{[\sup(g^\prime)+\epsilon \lambda_i\inf(f^\prime)]t}$. Thus, we have the following result.
	\begin{prop}
		Let $\dc$ be a discrete dynamical network with the hyperedge weight and the vertex measure of $H$ are $w_H$, $m_H$, respectively, and $f,g$ are differentiable functions with bounded derivatives,  $\sup(f^\prime)=\sup\limits_{x\in\mathbb{R}}(\frac{df(x)}{dx})$, $\sup(g^\prime)=\sup\limits_{x\in\mathbb{R}}(\frac{dg(x)}{dx})$. If all the eigenvalues of $ C_{(H,m_H,w_H)}$ belongs to the interval $(\infty, -\frac{\sup(g^\prime)}{\epsilon \inf(f^\prime)}]$ then synchronization in $\dc$ is stable under small perturbations.
	\end{prop}
	\section{Numerical illustrations}\label{sec-4}
	In this section, we numerically demonstrate the theoretical results obtained in the previous sections.
	\subsection{Comparison with graph-models}\label{s-ns}
\begin{figure}[ht]
		\centering
		\begin{subfigure}[b]{0.4\linewidth}
			\begin{tikzpicture}
				\filldraw[](0,-2)circle (2.5pt);
				\draw[](0,-1.5) node{$v_1$};
				\filldraw[](1,0)circle (2.5pt);
				\draw[](1,0.5) node{$v_2$};
				\filldraw[](2,-2)circle (2.5pt);
				\draw[](2,-1.5) node{$v_3$};
				\draw[](1,-1.2)circle (55pt);
				\draw[dotted, very thick](0,-2)--(1,0)--(2,-2)--cycle;
			\end{tikzpicture}
			\caption{A $3$-hyperedge is approximated by a triangle graph.}
			\label{fig:triangle}
		\end{subfigure}
		\begin{subfigure}[b]{0.4\linewidth}
			\begin{tikzpicture}[scale=1.3]
				\filldraw[](0,-1.5)circle (2.5pt);
				\draw[](-0.3,-1.8) node{$v_1$};
				\filldraw[](1,-2)circle (2.5pt);
				\draw[](0.6,-2) node{$v_2$};
				\filldraw[](2,-2)circle (2.5pt);
				\draw[](2.4,-2) node{$v_3$};
				\filldraw[](3,-1.5)circle (2.5pt);
				\draw[](3,-1.7)node{$v_4$};
				\filldraw[](0.5,-3)circle (2.5pt);
				\draw[](0.8,-3) node{$v_5$};
				\filldraw[](1,-4)circle (2.5pt);
				\draw[](0.5,-4) node{$v_6$};
				\filldraw[](2.5,-3)circle (2.5pt);
				\draw[](2.25,-3) node{$v_7$};
				\filldraw[](2,-4)circle (2.5pt);
				\draw[](2.5,-4) node{$v_8$};
				\draw [rounded corners=20](-0.5,-2.5) rectangle (3.5,-1);
				\draw [rounded corners=10](1.2,-1.3) rectangle (0.3,-4.5);
				\draw [rounded corners=10](1.7,-1.3) rectangle (2.7,-4.5);
				\draw[dashed, very thick](1,-2)--(1,-4);
				\draw[dashed, very thick](0,-1.5)--(3,-1.5);
				\draw[dashed, very thick](2,-2)--(2,-4);
				\draw[dashed, very thick](0,-1.5)--(1,-2)--(2,-2)--cycle;
				\draw[dashed, very thick](1,-2)--(3,-1.5)--(2,-2);
				
				\draw[dashed, very thick](1,-2)--(0.5,-3)--(1,-4);
				\draw[dashed, very thick](2,-2)--(2.5,-3)--(2,-4);
			\end{tikzpicture}
			\caption{Approximation of a hypergraph topology by a graph}
			\label{fig:hyp}
		\end{subfigure}
		\caption{Comparison between a hypergraph and its underlying graph}
	\end{figure}
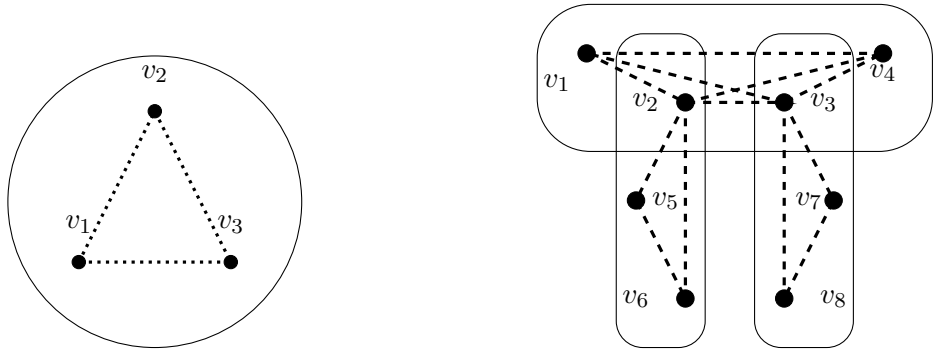

	We intend to incorporate multi-body interactions using the underlying topology of the dynamical network in our work. Generally, a network, that is, a graph, is used as the underlying topology of dynamical networks. Here we use hypergraph as the underlying topology. A hypergraph is a generalization of a graph, and if we take $cr(H)=rk(H)=2$, it becomes a graph. %Thus our model is so adaptable that we can return to the graph case when needed. 
	It is important to note that in many situations, the interactions in a network are multi-nary. Using binary interactions, we approximate those multi-nary synergies. For example, the interactions are multi-nary in the synchronized chirping of crickets and the synchronous flashing of a swarm of male fireflies. The diffusion of any substance (e.g., ink) on a surface ( a piece of cloth) can be described conveniently using grids on that surface. The grids are indeed hypergraphs, in which each cell is a vertex, and each cell, along with all its neighbouring cells, forms a hyperedge. Suppose  we put a drop of ink on one of the cells in the grid. As the ink spreads in all the neighbouring cells, the interaction is multi-nary and can be described conveniently by a hyperedge containing the cell and its neighbours. Though sometimes, multi-nary interaction can be approximated by multiple binary interactions, this approximation may not work in some situations.
	To illustrate this, we consider an abstract example where a $3$-uniform hyperedge, $e$, is approximated by three $2$-edges (drawn with the dotted lines in \Cref{fig:triangle}). The corresponding diffusion matrix  $L_G$ (which is the negative of Laplacian) of the triangle is $\left(\begin{smallmatrix}-2 &\phantom{-}1&\phantom{-}1 \\ \phantom{-} 1&-2&\phantom{-}1\\\phantom{-}1&\phantom{-}1 &-2\end{smallmatrix}\right)$.  The matrix representation of the operator $C_{(e,m_e,w_e)}$ with $m_e(v)=1$ for all $v$, and $w_e(e)=1$, corresponding to the edge $e$ is $\left(\begin{smallmatrix}     -1 &\phantom{-}\frac{1}{2}&\phantom{-}\frac{1}{2} \\ \phantom{-} \frac{1}{2} &-1 &\phantom{-} \frac{1}{2} \\ \phantom{-} \frac{1}{2}&\phantom{-}\frac{1}{2}  &-1 \end{smallmatrix}\right)$. Now  the discrete difference equation expressed by $L_G$ is \begin{align}\label{disL}
		x_{t+1}=x_t+ L_G(x_t),
	\end{align}
	and the same described by $C_{(e,m_e,w_e)}$ is
	\begin{align}
		\label{disH}
		x_{t+1}=x_t+C_{(e,m_e,w_e)}(x(n)).
	\end{align}
	The eigenvalues of $I_3+L_G$ are $-2,-2,1$ and the same of $I_3+C_{(e,m_e,w_e)}$ are $-\frac{1}{2},-\frac{1}{2},1$. 	 The trajectories of the system given by \Cref{disH} synchronize because for the matrix $I_3+C_{(e,m_e,w_e)}$  the eigenspace of $1$ is the vector space generated by the vector $\mathbf{1}$ and the  magnitude of the other eigenvalues is less than $1$. { The same is not valid for \Cref{disL} because the magnitudes of the eigenvalues 
		(other than $1$) of $I_3+L_G$ are greater than $1$. }
	
	% The matrix $C_{(e,m_e,w_e)}$ corresponding to the hypergraph in  \Cref{fig:triangle} coincides with the negative normalized Laplacian matrix (defined in \cite[Equation-2]{1}) of the triangle graph. In general, the matrix $C_{(e,m_e,w_e)}$ is different from this notion of normalized Laplacian because $C_{(e,m_e,w_e)}$ is symmetric for any hypergraph. In contrast, the normalized Laplacian is not symmetric for non-regular graphs.
	Next, we consider another abstract example where binary interactions fail to approximate multi-nary interactions. We approximate the hypergraph by a graph whose edges are drawn with the dashed lines (see \Cref{fig:hyp}).
	The diffusion operator $C_{(H,m_H,w_H)}$ of the hypergraph and the negative Laplacian matrix $L_G$ of the underlying graph $G$ are as follows.
	$C_{(H,m_H,w_H)}=\left(\begin{smallmatrix}
		-1 & \frac{1}{3} & \frac{1}{3} & \frac{1}{3} & 0 & 0 & 0 & 0\\ \frac{1}{3} & -2 & \frac{1}{3} & \frac{1}{3} & \frac{1}{2} & \frac{1}{2} & 0 & 0\\ \frac{1}{3} & \frac{1}{3} & -2 & \frac{1}{3} & 0 & 0 & \frac{1}{2} & \frac{1}{2}\\ \frac{1}{3} & \frac{1}{3} & \frac{1}{3} & -1 & 0 & 0 & 0 & 0\\ 0 & \frac{1}{2} & 0 & 0 & -1 & \frac{1}{2} & 0 & 0\\ 0 & \frac{1}{2} & 0 & 0 & \frac{1}{2} & -1 & 0 & 0\\ 0 & 0 & \frac{1}{2} & 0 & 0 & 0 & -1 & \frac{1}{2}\\ 0 & 0 & \frac{1}{2} & 0 & 0 & 0 & \frac{1}{2} & -1 
	\end{smallmatrix}\right)$
	, $L_G=\left(\begin{smallmatrix}
		-3 & 1 & 1 & 1 & 0 & 0 & 0 & 0\\ 1 & -5 & 1 & 1 & 1 & 1 & 0 & 0\\ 1 & -5 & 1 & 1 & 0 & 0 & 1 & 1\\ 1 & 1 & 1 & -3 & 0 & 0 & 0 & 0\\ 0 & 1 & 0 & 0 & -2 & 1 & 0 & 0\\ 0 & 1 & 0 & 0 & 1 & -2 & 0 & 0\\ 0 & 0 & 1 & 0 & 0 & 0 & -2 & 1\\ 0 & 0 & 1 & 0 & 0 & 0 & 1 & -2 
	\end{smallmatrix}\right)
	$.

	Now, we compare the diffusion equations stated using $C_{(H,m_H,w_H)}$ and $L_G$, respectively.
	\begin{align}
		\label{disC2}
		x_{t+1}=x_t+\frac{3}{4}C_{(H,m_H,w_H)}(x_t).
	\end{align}
	\begin{align}\label{disL2}
		x_{t+1}=x_t+\frac{3}{4}L_G(x_t).
	\end{align}
	The eigenvalues of the matrix, $I_8+\frac{3}{4}C_{(H,m_H,w_H)}$ are $  -0.930$%77774934061352496428298763931\\
	, $-0.678$%05361261225594038393182927393\\  
	, $-0.125$,%\\
	$ -0.125$, $0.867\times 10^{-17}$%\\ , 0.0000000000000000086736173798840354720596224069595\\ 
	, $0.553$,%0.55305361261225660651774660436786\\ 
	$0.806$, %0.80577774934061296985277067506104\\ 1.0 >
	$1$. Here $\mathbf{1}$ is the eigenvector corresponding to the eigenvalue $1$, and the  absolute values of the other eigenvalues are less than $1$. Therefore, the trajectories of the system given by \Cref{disC2} synchronize asymptotically.
	The eigenvalues of the matrix, $I_8+\frac{3}{4}L_G$ are
	$-3.17$, $ -3.79$, $-2 $, $-1.25 $, $-1.25$, $-0.08$, $0.53$, $1 $. So the absolute values of some eigenvalues are greater than $1$. Therefore, the trajectories of the system given by \Cref{disL2} may not synchronize.
	\subsection{Chemical-gene interaction (biogrid) and protein complex hypergraphs }\label{comapare}
	Now we extend our study to two hypergraphs, (i) bio grid hypergraph and (ii) protein-complex hypergraph, created from real data. 
	Biogrid hypergraph is constructed from human chemical-gene(target) interactions which are useful for studying drug-gene(target) interactions (The data is downloaded on 18/02/2019 from the  repository, BioGRID \cite{chatr2017biogrid}). Here the genes are considered as vertices, and chemicals are as hyperedges. A hyperedge corresponding to a chemical is constituted by a group of genes that are targeted by that chemical.
	A protein complex hypergraph is created from the database, CORUM  \cite{giurgiu2019corum}, a resource of mammalian protein complexes (the data is also downloaded on 18/02/2019) and which is also useful for predicting unknown interactions between proteins. Here protein complexes are considered as vertices and subunits are as hyperedges. A hyperedge(subunit) is constructed with the protein complex associated with the corresponding subunit. 
	Initially, the biogrid hypergraph contained $2138$ vertices and $4455$ hyperedges, whereas the protein complex hypergraph was made of $3638$ vertices and $ 2848$ hyperedges. 
	After removing all the hyperedges containing only one vertex, we find $1501$ hyperedges in our biogrid hypergraph. Since our theoretical results are on connected hypergraph, we use the largest connected component %(retrieved by using \url{https://www.mathworks.com/matlabcentral/fileexchange/30926-largest-component})
	as the underlying topology of  dynamical networks in our study. The largest connected component of the biogrid hypergraph consists of $1808$ vertices and $ 1431$ hyperedges. The same of the protein complex hypergraph contains   $2770$ vertices and $2383$ hyperedges. 
	%%%%%%%%%%%
	\subsubsection{Comparison with the hypergraph and its underlying graph} 
	Now we consider the diffusion equations,
	\begin{align}
		\label{disCbio}
		x_{t+1}=x_t+\frac{1}{110}C_{(H,m_H,w_H)}(x_t)
	\end{align}
	involving the diffusion operator $C_{(H,m_H,w_H)}$, and 
	\begin{align}\label{disLbio}
		x_{t+1}=x_t+\frac{1}{110}L_G(x_t)
	\end{align}
	containing the negative Laplacian of the underlying graph of the biogrid hypergraph.
	\begin{figure}[H]
		\begin{subfigure}{0.45\linewidth}
			\includegraphics[scale=0.4]{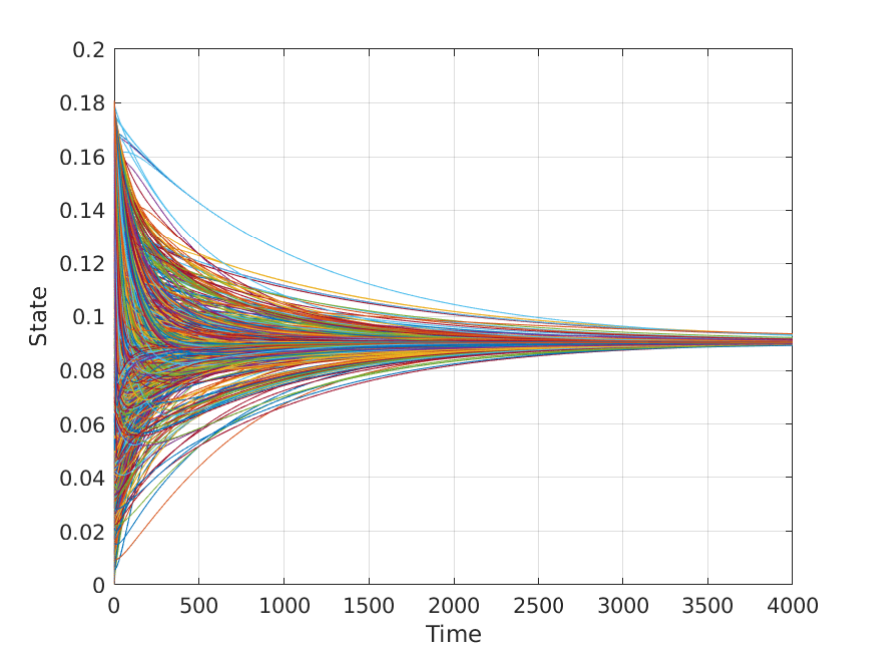}
			\caption{Trajectories of the dynamical system given by \Cref{disCbio}.}
			\label{fig:biogrid_C}
		\end{subfigure}
		\begin{subfigure}{0.45\linewidth}
			\includegraphics[scale=0.4]{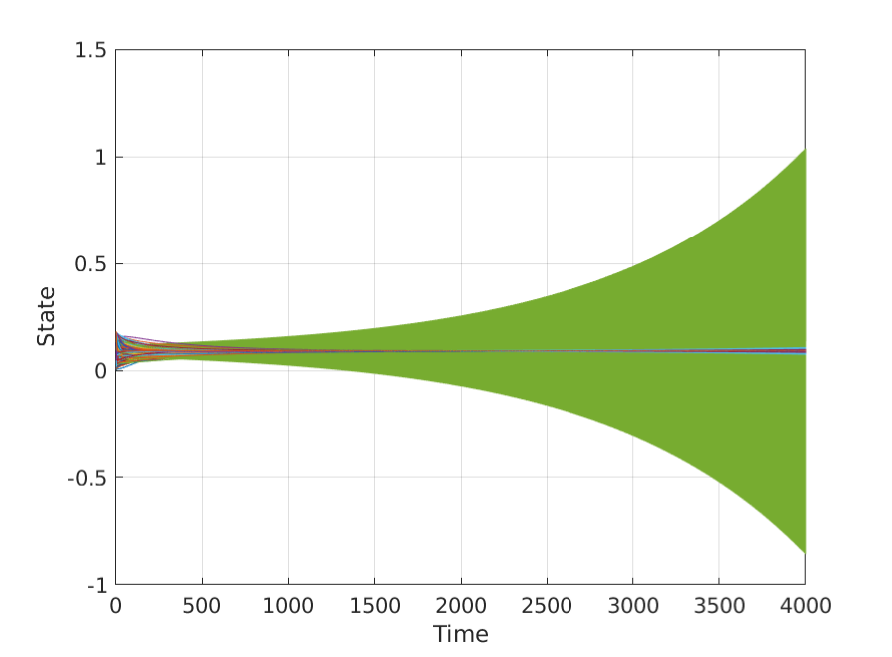}
			\caption{Trajectories of the dynamical system given by \Cref{disLbio}.}
			\label{fig:biogrid_L}
		\end{subfigure}
		\caption{ Comparison with the hypergraph and its underlying graph in Chemical-gene interaction (biogrid) hypergraph. }
	\end{figure}
	
\begin{figure}[H]
		\begin{subfigure}{0.45\linewidth}
			\includegraphics[scale=0.4]{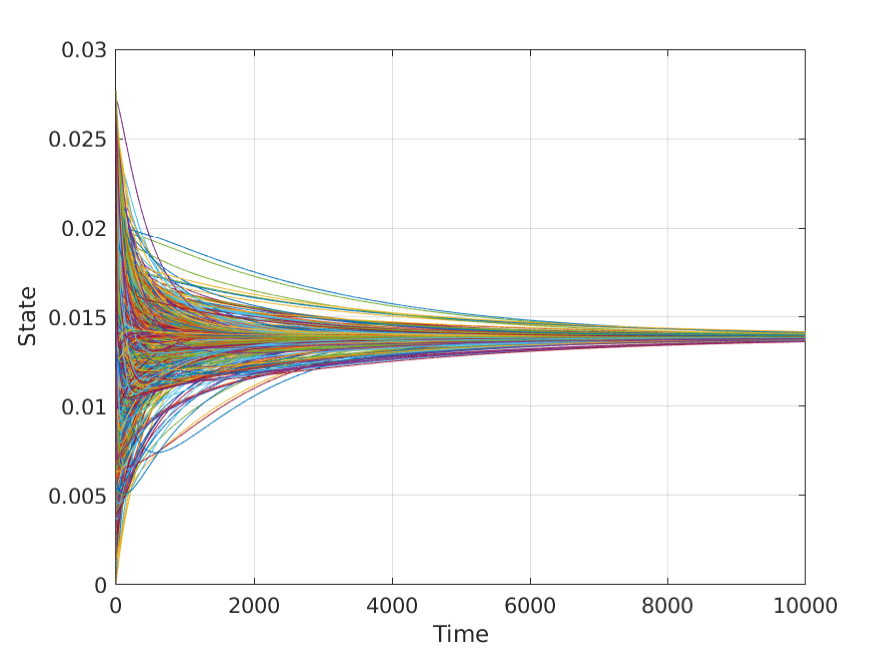}
			\caption{Trajectories of the dynamical system given by \Cref{disCprot}}
			\label{fig:prot_C}
		\end{subfigure}
		\begin{subfigure}{0.45\linewidth}
			\includegraphics[scale=0.4]{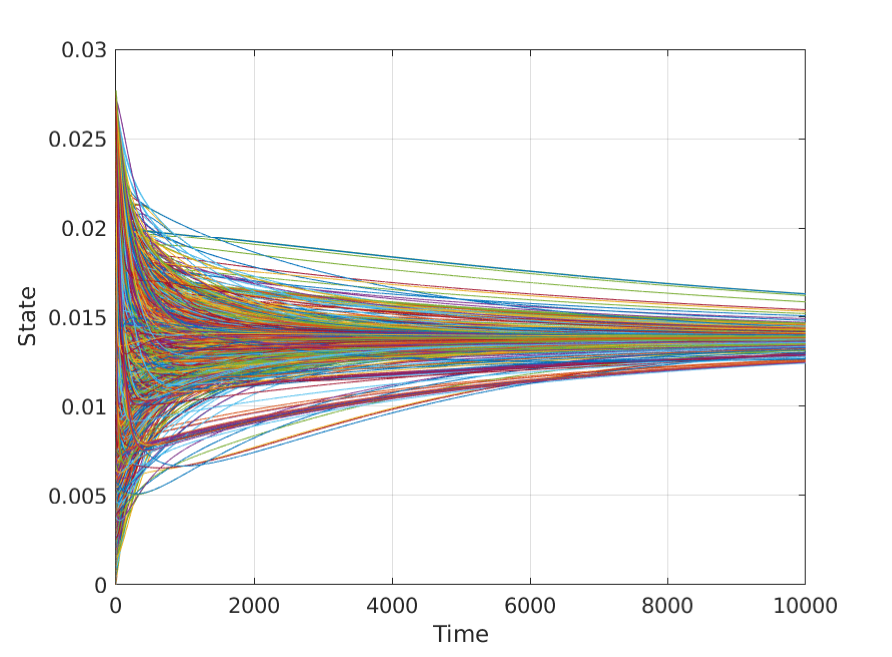}
			\caption{Trajectories of the dynamical system given by \Cref{disLprot}.}
			\label{fig:prot_L}
		\end{subfigure}
		\caption{Comparison with the hypergraph and its underlying graph in protein complex hypergraph.   }
	\end{figure}

	The converging evolution of the trajectories (\cref{fig:biogrid_C}) given by \Cref{disCbio} demonstrates the presence of a diffusion process in the dynamical system. In contrast, the evolution of the trajectories given by \Cref{disLbio} is divergent in nature (\cref{fig:biogrid_L}). Therefore, for biogrid hypergraph, $C_{(H,m_H,w_H)}$ is a better diffusion operator than $L_G$.
	To have another instance where $C_{(H,m_H,w_H)}$ is proved to be a better diffusion operator than the negative laplacian $L_G$ of the underlying graph of the same hypergraph, we present one more comparison.
	Consider the following discrete difference equations
	\begin{align}
		\label{disCprot}
		x_{t+1}=x_t+\frac{1}{200}C_{(H,m_H,w_H)}(x_t)
	\end{align}
	involving the diffusion operator $C_{(H,m_H,w_H)}$ of the protein complex hypergraph and 
	\begin{align}\label{disLprot}
		x_{t+1}=x_t+\frac{1}{200}L_G(x_t)
	\end{align}
	involving the negative Laplacian of the underlying graph of the protein complex hypergraph. Considering the $10000$ iterations in the trajectories of the dynamical system given by \Cref{disCprot} (\Cref{fig:prot_C}) and \Cref{disLprot}  (\Cref{fig:prot_L}), it is clear that the trajectories of the system given by \Cref{disCprot} converging faster than that of the system given by \Cref{disLprot}.
	
	\subsubsection{Global analysis of synchronization in dynamical networks with biogrid hypergraph topology}	In this section we will verify our theoretical results on global synchronization with some dynamical networks with the biogrid hypergraph as its underlying architecture. In the following examples we consider $C_{(H,m_H,w_H)}$ with $m_H(v)=1$ for all $v\in V(H)$, and $w_H(e)=1$ for all $e\in E(H)$.
	\begin{exm}
		In this example, we consider the dynamical network with the biogrid hypergraph as its underlying topology. 
		
		a)We set $k=1  $. If we define $\bar{f}:\mathbb{R}\to\mathbb{R}$  as $x\mapsto q\sin{(-x)}$ and  $\bar{g}:\mathbb{R}\to\mathbb{R}$  as $x\mapsto p\cos{(-x)}$ then the lipschitz constants $k_g=p,k_f=q$. If we choose $\epsilon=\frac{1}{88}$ and $p=0.4,q=0.5$ then $  [k_g+\epsilon \|C_{(H,m_H,w_H)}\|k_f]=(p+\frac{1}{88}(87.6182)q)< 1$, which is the condition given in \Cref{lpshz-dis}. Thus, the trajectories of the dynamical network synchronize (\Cref{fig:bio_sync}). If we choose $p=1,q=1.53$, then the condition given in \Cref{lpshz-dis} does not satisfied and the trajectories remain asynchronous (\Cref{fig:bio_async}). If $p=1,q=1.52$ then $  [k_g+\epsilon \|C_{(H,m_H,w_H)}\|k_f]>1$ that is the conditions  of \Cref{lpshz-dis} are not satisfied. Despite that, the trajectories synchronize (\Cref{fig:bio_sync_ne}). This shows that the condition is sufficient but not necessary.
		
		\begin{figure}[htbp]
			\centering
			\begin{subfigure}{0.3\textwidth}
				\includegraphics[scale=0.4]{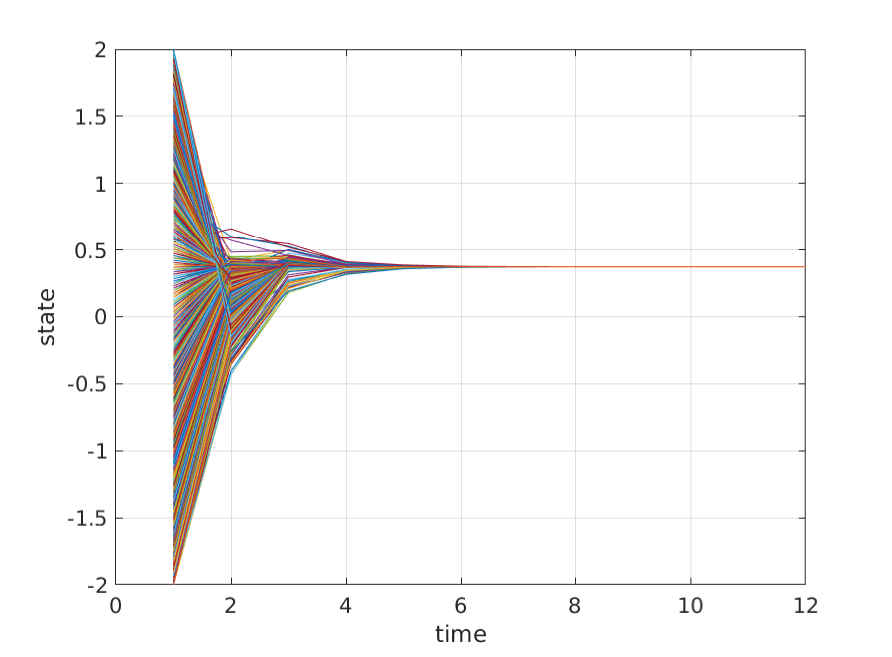}
				\subcaption{ Synchronization with $  [k_g+\epsilon \|C_{(H,m_H,w_H)}\|k_f]<1$.}
				\label{fig:bio_sync}
			\end{subfigure}\hfill
			\begin{subfigure}{0.3\textwidth}
				\includegraphics[scale=0.4]{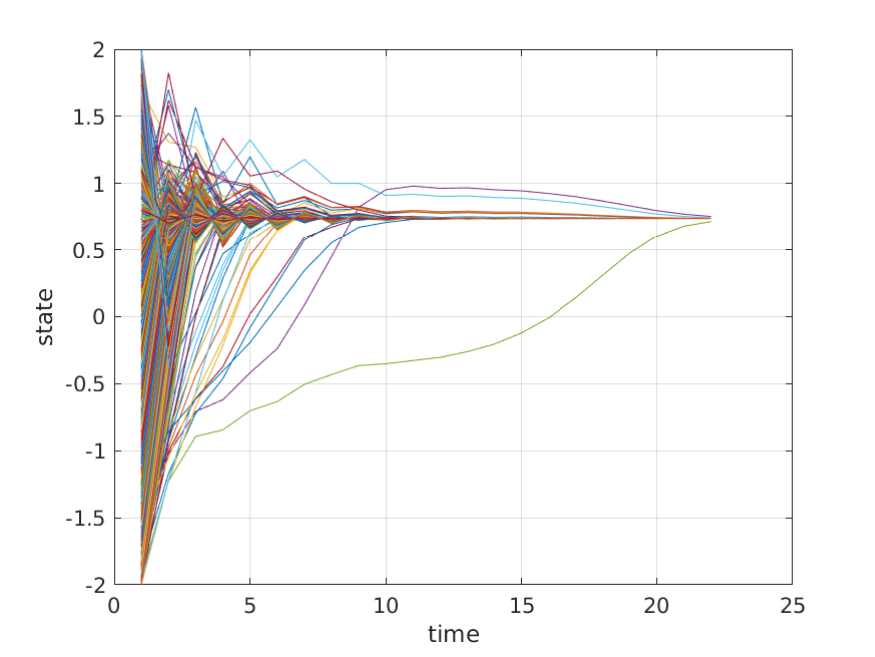}
				\subcaption{ Synchronization with $  [k_g+\epsilon \|C_{(H,m_H,w_H)}\|k_f]>1$.}
				\label{fig:bio_sync_ne}
			\end{subfigure}\hfill
			\begin{subfigure}{0.3\textwidth}
				\includegraphics[scale=0.4]{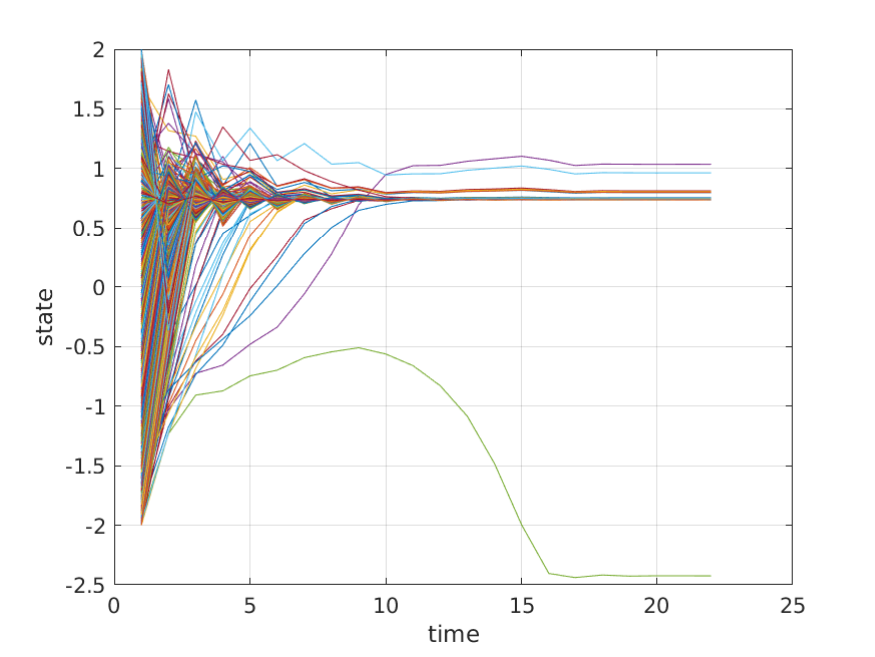}
				\subcaption{Asynchronous trajectories with $  [k_g+\epsilon \|C_{(H,m_H,w_H)}\|k_f]>1$.}
				\label{fig:bio_async}
			\end{subfigure}
			\caption{Simulation of \Cref{lpshz-dis}.}
		\end{figure}
		%\end{exm}\begin{exm}In this example, we consider the dynamical network with the biogrid hypergraph as its underlying topology.
		
		b)If the coupling strength $\epsilon=\frac{1}{45}$, $k=1  $ then  $\|[I_{|V|}+\epsilon C_{(H,m_H,w_H)}]\|=1$. Now if we define $\bar{f}=\bar{g}:\mathbb{R}\to\mathbb{R}$  as $x\mapsto qe^{\sin{x}}$, then $k_f=k_g=\sup{\|\frac{d}{dx}(\bar{f}(x))}\|=\sup{\|(q\cos{(x)}e^{\sin{(x)}})}\|\le qe$. Therefore, if we choose $q=\frac{1}{2.8}<\frac{1}{e}$ then  $\|[I_{|V|}+\epsilon C_{(H,m_H,w_H)}]\|< \frac{1}{k_f}$ and $\|[I_{|V|}+\epsilon C_{(H,m_H,w_H)}]\|<\frac{1}{\|\sup\bar{f}^\prime\|}$ which are the condition given in \Cref{lpshz-dis}. Therefore, the trajectories synchronize (\Cref{fig:q28}). Calculating in Matlab we get $k_f=k_g=\sup{\|\frac{d}{dx}(\bar{f}(x))}\|=\sup{\|(q\cos{(x)}e^{\sin{(x)}})}\|\approx q\times1.46$. Thus, if we choose $q=\frac{1}{1.47}$, it agrees with the condition given in given in \Cref{lpshz-dis}. Therefore, the trajectories synchronize (\Cref{fig:q147}). When $q=\frac{1}{1.15}$, the condition given in \Cref{lpshz-dis} is not satisfied and the trajectories remain asynchronous (\Cref{fig:q115}). However, if $q=\frac{1}{1.15}$ then also the condition is not satisfied but synchronization is observed in this case (\Cref{fig:q12}). This is because the condition is sufficient but not necessary.
		
		\begin{figure}[htbp]
			\begin{subfigure}{.45\textwidth}
				\centering
				% include first image
				\includegraphics[scale=0.25]{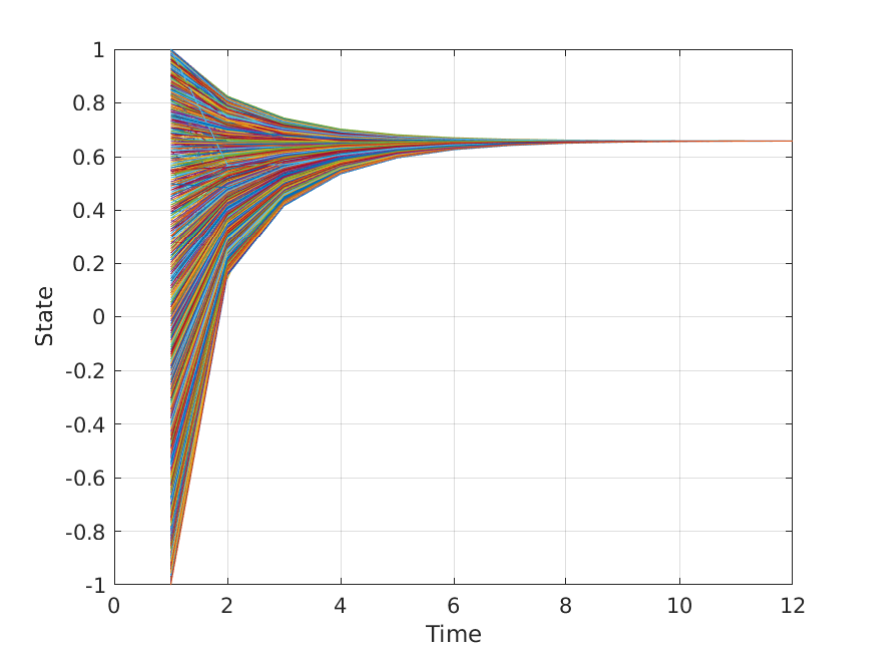}  
				\caption{Synchronization with $q=\frac{1}{2.8}$. }
				\label{fig:q28}
			\end{subfigure}
			\begin{subfigure}{.45\textwidth}
				\centering
				% include first image
				\includegraphics[scale=0.25]{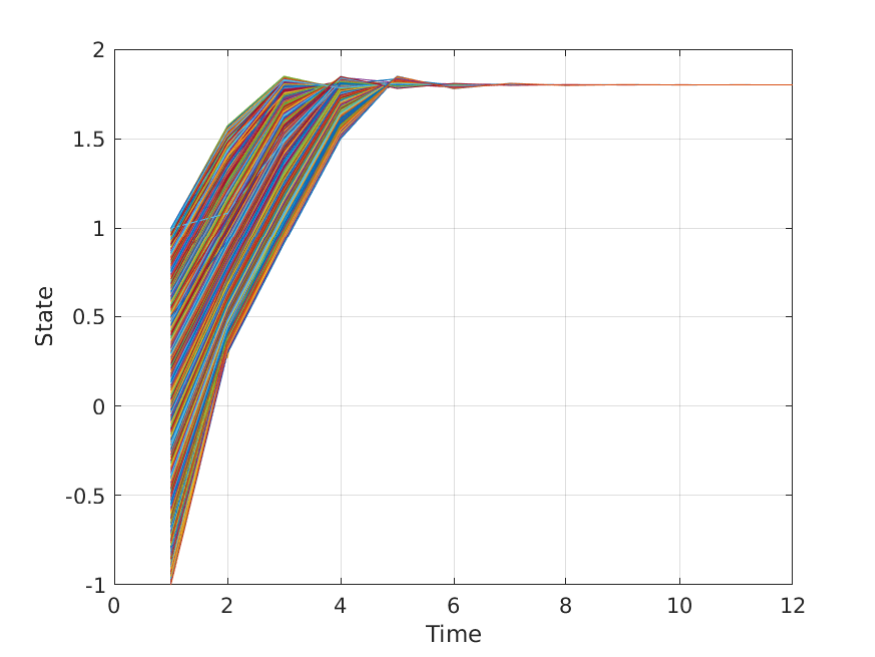}  
				\caption{Synchronization with $q=\frac{1}{1.47}$.}
				\label{fig:q147}
			\end{subfigure}
			\begin{subfigure}{.45\textwidth}
				\centering
				% include second image
				\includegraphics[scale=0.25]{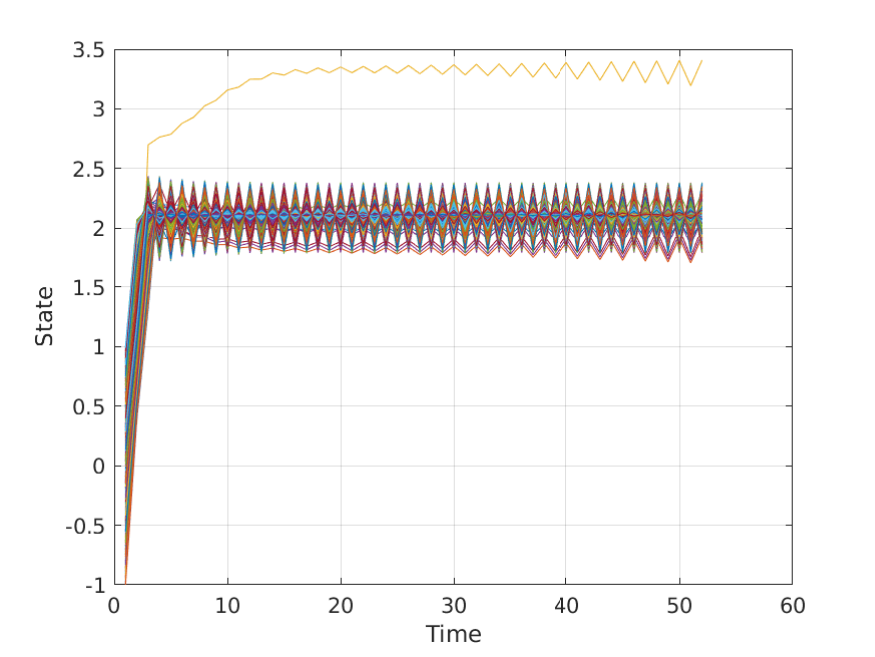}  
				\caption{Asynchronous trajectories with $q=\frac{1}{1.15}$}
				\label{fig:q115}
			\end{subfigure}
			\begin{subfigure}{.45\textwidth}
				\centering
				% include second image
				\includegraphics[scale=0.25]{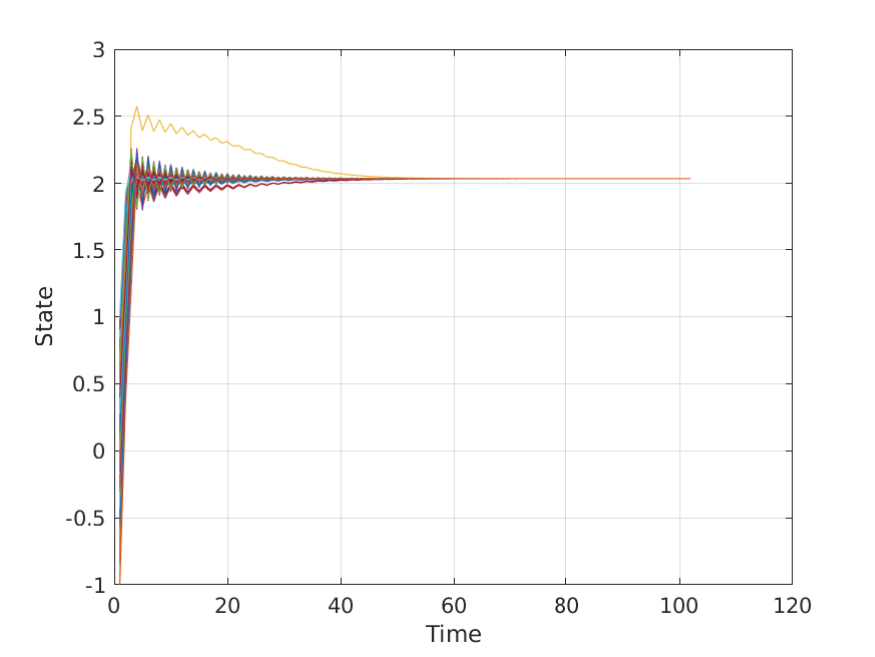}  
				\caption{Synchronization with $q=\frac{1}{1.2}$}
				\label{fig:q12}
			\end{subfigure}
			\caption{Simulation of \Cref{lpshz-dis}.}
			\label{fig:fig}
		\end{figure}
	\end{exm}

	\subsection{Approximation of a hypergraph by a weighted graph}
	The previous examples show that the diffusive actions of a hypergraph and its underlying projection graph are not the same. Now we show how the underlying hypergraph representation can be well approximated by a weighted  graph for studying diffusion. Let the negative Laplacian of a weighted graph $G$,  with edge-weight $w_G:E(G)\to (0,\infty)$. The negative Laplacian $L_{(G,m_G,w_G)}:\mathbb{R}^{V(G)}\to \mathbb{R}^{V(G)}$ of $G$ is defined as, for all $x\in \mathbb{R}^{V(G)}$, and $v\in V(H)$, $(L_{(G,m_G,w_G)}x)(v)=\sum\limits_{v\in V(G)}\frac{a_{uv}}{m_H(u)}(x(v)-x(u))$, where $a_{uv}=w_H(\{u,v\})$ if $\{u,v\}\in E(G)$, otherwise $a_{uv}=0$. 
	\begin{prop}
		Let $H$ be a hypergraph with vertex measure $m_H$ and hyperedge weight $ w_H$. If $G $ is the underlying projected graph of $H$ with $m_G=m_H$, and $w_G$ is defined as for all $ \{u,v\}\in E(H)$, $w_G(\{u,v\})=\sum\limits_{u\in V(H)}\sum\limits_{e\in E_u(H)\cap E_v(H)}\frac{w_H(e)}{|e|-1}$, then $C_{(H,m_H,w_H)}=L_{(G,m_G,w_G)}$.
	\end{prop}
	\begin{proof}
		For all $x\in \mathbb{R}^{V(H)} $, and $v\in V(H)$,
		\begin{align*}
			( C_{(H,m_H,w_H)}(x))(v)&=\sum\limits_{i=cr(H)}^{rk(H)}\frac{1}{i-1}\sum\limits_{e\in E_v(H_i)}\frac{w_{H_i}(e)}{m_{H_i}(v)}\sum\limits_{u\in e}(x(u)-x(v))\\
			&=\sum\limits_{e\in E_v(H)}\frac{w_{H}(e)}{m_{H}(v)}\frac{1}{|e|-1}\sum\limits_{u\in e}(x(u)-x(v))\\
			&=\frac{1}{m_H(v)}\sum\limits_{u\in V(H)}\left(\sum\limits_{e\in E_u(H)\cap E_v(H)}\frac{w_H(e)}{|e|-1}\right)(x(u)-x(v))\\
			&=(L_{(G,m_G,w_G)}(x))(v).
		\end{align*}  
		Thus, the result follows.
	\end{proof}
	The action of the diffusion operator is the same as a weighted version of the underlying projected graph. The information of the hypergraph is encoded in the weight of the graph. Therefore, the unweighted version of the underlying projected graph is not a good approximation for a hypergraph considering an underlying structure of a dynamical network for studying diffusion. 
	
	\section{Discussion }	In this paper, our focus was synchronization in a dynamical network. To incorporate the multi-body interaction in a dynamical network, we have used here a hypergraph $H$ with vertex measure $m_H$ and hyperedge weight $w_H$.	we introduced the operator $C_{(H,m_H,w_H)}$ associated with $H$, which acts as a diffusion operator in our dynamical networks model. For the same hypergraph $H$, for different choices of $m_H$, $w_H$, the operator $C_{(H,m_H,w_H)}$ provide us different diffusion operators associated with $H$.
	
	For example, let us consider the dynamical network of heat propagation, where the total amount of heat in the network is constant. We use $\delta_{V(H)}$ as a normalizing factor.
	As another example, to get the negative of the normalized Laplacian described in \cite{MR4208993}, we choose $m_{H}(v)=|E_v(H)|$. On the other hand, we consider $m_{H}$ as a constant function in the dynamical network of malware propagation on the internet, where malware replicates itself. Thus the total amount of malware in the network grows with time.
	We set $w_{H}(e)=|e|-1$ for all $e \in E(H)$ when the cardinality of the hyperedges in the hyperedge-coupling is not relevant in a dynamical network. In this case, our diffusion operator becomes the negative of the hypergraph Laplacian considered in \cite{rodriguez2003Laplacian,rodriguez2009Laplacian,bretto2013hypergraph}. 
	%Thus, here $\sigma_H(e)=1$ for all $e\in E(H)$ and the difference of hyperedge cardinality does not affect the diffusion operator. 
	Whereas  in \cite{MR4208993}, the value of $m_H(e)$ is taken as $1$. The negative multiple of the Laplacian operator. 
	% \subsection{Global analysis of synchronization in continuous dynamical networks}
	
	There are two main types of diffusion phenomenon on networks. The first one is when there are several dynamical systems on the vertices (nodes) of the network, and two connected nodes  affect the dynamics of each other through the hyper(edges), which act as diffusion coupling. We have already described this phenomenon, and while doing so, we have developed the diffusion operator $C_{(H,m_H,w_H)}$. One can use this operator to  explain another type of diffusion phenomenon in networks. In this phenomenon, there is no dynamical network on the nodes, but each node contains some substance (e.g.fluid, information, disease), and the substance diffuses through the hyperedges. This second type of phenomenon is called \emph{random walk}. A random walk on hypergraph $H$ is a map $\mathfrak{r}_H:\mathbf{T}\to V(H)$, such that $\mathbf{T}=\mathbb{N}\cup\{0\}$,%set of all non-negative integers, is the domain of time
	and $\mathfrak{r}_H(t)$ depends only on its previous state $\mathfrak{r}_H(t-1)$ for all $t(\ne 0)$. Thus, the event $\mathfrak{r}_H(uv)=(\mathfrak{r}_H(t+1)=v|\mathfrak{r}_H(t)=u)$ is independent of $t$. For each $e\in E_u(H)\cap E_v(H)$, we denote the event of the random walk $\mathfrak{r}_H$ going from $u$ to $v$ through the hyperedge $e$ as $\mathfrak{r}_H(uev)$. Thus, $\mathfrak{r}_H(uv)=\sum\limits_{e\in E_u(H)\cap E_v(H)} \mathfrak{r}_H(uev)$. Again, $\mathfrak{r}_H(uev)=\mathfrak{r}_H(ue)\mathfrak{r}_H(ev)$. Therefore, $Prob(\mathfrak{r}_H(uev))=Prob(\mathfrak{r}_H(ue))Prob(\mathfrak{r}_H(ev)|\mathfrak{r}_H(ue))$. Since $\mathfrak{r}_H(ue)$ depends on $u$ and $e$, we can define $m_H$, and $w_H$ in such a way that $Prob(\mathfrak{r}_H(ue))=\frac{w_H(e)}{m_H(u)}$. Since in $(\mathfrak{r}_H(ev)|\mathfrak{r}_H(ue))$, the random walker is coming from $u$ to $e$, and random walker can not stay in a particular vertex for two subsequent time steps (that is the random walk is non-lazy), form $e$, it has to choose one of the remaining $|e|-1$ vertices. Thus, $Prob(\mathfrak{r}_H(ev)|\mathfrak{r}_H(ue))=\frac{1}{|e|-1}$. Therefore, $Prob(\mathfrak{r}_H(uv))=\sum\limits_{e\in E_u(H)\cap E_v(H)}\frac{w_H(e)}{m_H(u)}\frac{1}{|e|-1}$. Since $\bigcup\limits_{v\in V(H)}\mathfrak{r}_H(uv)$ is a certain event, $Prob(\bigcup\limits_{v\in V(H)}\mathfrak{r}_H(uv))=1$. Thus, $m_H(u)=\sum\limits_{e\in E_u(H)}\frac{w_H(e)}{|e|-1}$, for all $u\in V(H)$. The probability transition matrix $P=(P_{uv})_{u,v\in V(H)}$ is defined as $P_{uv}=Prob(\mathfrak{r}_H(uv))$ for $u\ne v$, and otherwise $P_{vv}=0$. Since $m_H(u)=\sum\limits_{e\in E_u(H)}\frac{w_H(e)}{|e|-1}$, for all $x\in\mathbb{R}^{V(H)}$, $C_{(H,m_H,w_H)}x=Px-x$.
	
	%\section{Application}
	%Now we use our models and results to explore some natural phenomena.
	%\subsection{Ascent of sap}
	
	The upward transport of water and minerals from the root to an upper section of a plant body using the xylem and  phloem tissue is known as the ascent of sap in plants. The xylem is a complex tissue of both living and nonliving cells, and the phloem is a complex living tissue. Suppose we consider cells of the xylem and phloem tissue to be the vertices (nodes) of a hypergraph. Each node (cell), along with all its neighbouring cells, forms a hyperedge. These hyperedges  form a connected system of water-conducting channels reaching all parts of the plants. In this case, for any time $n\in \n$, $x(n)(u)$ is a (multi-dimensional)
	variable which keeps track of various information, such as \textit{water potential}, the density of sap, sucrose concentration in the $u$-th cell at time $n$. Water potential measures the tendency of water to migrate from one location to another owing to osmosis, gravity, and other factors. The difference of $x(n)$ in adjacent vertices creates a diffusion gradient (osmotic gradient) that draws water into the vertex from the adjacent vertices with higher values of $x(n)$. The function ${f}={g}$ can incorporate all the factors affecting the changes of $x(n)$ in each cell, such as evaporation and absorption of water and minerals, etc. After absorbing water in the root, the value of $x(n)$ in the nodes (cells) present in the root becomes higher, and sap transportation begins. The process ends in a synchronized state, in which the concentration of sap becomes equal in each node.
	%\subsection{Global wind flow}

	Latitude and longitude divide the world into grids. Let us consider a hypergraph in which each grid cell is a vertex, and a vertex with all its neighbours forms a hyperedge. For any time $n\in\n$, $x(n)$ is the state of the dynamical network,  such that $x(n)(u)$ is the air pressure at the vertex $u$. The function ${f}$ incorporates all the factors that regulate the air pressure of a region, such as elevation or altitude, the  average temperature, air composition, amount of water vapour, etc. Here ${f}={g}$. Under the hypergraph's diffusive influence, air flows from high-pressure regions to low-air-pressure regions to attain synchronisation. 
	%\subsection{Thermal conduction in solids} 
	
	Our dynamical network model can also explain the conduction of heat in solids. We can divide any solids into $3$-dimensional grids. Each cell of the grid is considered a vertex. A node, along with all its neighbouring cells, forms a hyperedge. In this case, $\bar{g}$ is the zero function, and ${f}$ is the identity function. $x(n)(u)$ is the average temperature at a cell $u$. The heat conduction terminates at a synchronized state.
	%\subsection{Managing financial accounts}
	
	Our dynamical network model can be used in the risk and profit management of financial investments of companies or individuals. Various types of investments have different risk and profit possibilities. Usually, the higher the profit, the greater the risk. Therefore, to optimize the profit, we may seek an algorithm to diffuse the risk of an investment with a higher return. We can construct a connected hypergraph in which the financial investments are vertices and investments of different risk, and profit levels are coupled to form the hyperedges. Here ${f}$ and ${g}$ can incorporate the gain of interest or loss in individual investments. The diffusive influence of the network diffuses the surplus money from the higher risk-higher profit investment to lower risk-lower profit investments and thus manages the risk. By changing the coupling strengths and the weights of the hyperedges, one can regulate the risks and profits.

	\section*{Acknowledgements}
	The authors would like to show their gratitude to Shirshendu Chowdhury for his constant valuable advice and guidance in studying dynamical systems throughout this work. His constructive comments have enriched this study enormously. The authors are also grateful to Timoteo Carletti for helpful discussions on the terminologies, hyper-networks, and hypergraphs.
	
	The authors sincerely thank Meghna, Chiranjeet Ghosh, Sourav Makhal, Saikat Mondal, Buddhadev Chatterjee, Kausik Das, Arun Sardar for their precious responses in a google form, that help us to construct a toy model of dynamical network with hypergraph as its underlying architecture. The authors are also thankful to Mrinmay Biswas, Sugata Ghosh, Gargi Ghosh, Avishek Chatterjee, Shibananda Biswas, and Somnath Basu for fruitful discussions.
	
	The work of the second author is supported by University Grants Commission, India (Beneficiary Code/Flag: 	BININ00965055 A).  
	\section*{Conflict of interest}
	The authors declare that they have no conflict of interest.
	\bibliographystyle{siam}
		
%	\bibliography{syn}		
\end{document}